\documentclass{amsart}
\usepackage{amssymb,latexsym}
\numberwithin{equation}{section}
\newtheorem{theorem}{Theorem}[section]
\newtheorem{proposition}[theorem]{Proposition}
\newtheorem{corollary}[theorem]{Corollary}

\newtheorem{remark}[theorem]{Remark}
\newtheorem{lemma}[theorem]{Lemma}

\begin{document}

\pagenumbering{arabic}
\pagestyle{headings}
\def\sof{\hfill\rule{2mm}{2mm}}
\def\ls{\leq}
\def\gs{\geq}
\def\SS{\mathcal S}
\def\qq{{\bold q}}
\def\txx{{\frac1{2\sqrt{x}}}}
\def\mn{\mbox{-}}
\def\sss{\mathcal{S}(1\mn3\mn2)}
\def\ssn{{\in S_n(1\mn\3\mn2)}}

\title{ {\sc  continued fractions and generalized patterns}}
   
\author{Toufik Mansour}
\maketitle
\begin{center}{LABRI, Universit\'e Bordeaux I,\\ 
              351 cours de la Lib\'eration, 33405 Talence Cedex\\
		{\tt toufik@labri.fr} }
\end{center}
%
%===========================================================================
\section*{Abstract}
In \cite{BS} Babson and Steingrimsson introduced generalized permutation patterns that
allow the requirement that two adjacent letters in a pattern must be adjacent in the 
permutation. 

Let $f_{\tau;r}(n)$ be the number of $1\mn3\mn2$-avoiding permutations 
on $n$ letters that contain exactly $r$ occurrences of $\tau$, where 
$\tau$ a generalized pattern on $k$ letters. Let $F_{\tau;r}(x)$ 
and $F_\tau(x,y)$ be the generating functions defined by 
$F_{\tau;r}(x)=\sum_{n\geq0} f_{\tau;r}(n)x^n$ and 
$F_\tau(x,y)=\sum_{r\geq0}F_{\tau;r}(x)y^r$.
We find an explicit expression for $F_\tau(x,y)$ in the form of a 
continued fraction for where $\tau$ given as a generalized pattern; 
$\tau=12\mn3\mn\dots\mn k$, $\tau=21\mn3\mn\dots\mn k$, 
$\tau=123\dots k$, or $\tau=k\dots 321$. 
In particularly, we find $F_\tau(x,y)$ for any $\tau$ generalized pattern 
of length $3$. This allows us to express $F_{\tau;r}(x)$ via Chebyshev 
polynomials of the second kind, and continued fractions.
%==========================================================================
\section{Introduction}
Let $[p]=\{1,\dots,p\}$ denote a totally ordered alphabet 
on $p$ letters, and let $\pi=(\pi_1,\dots,\pi_m)\in [p_1]^m$,
$\beta=(\beta_1,\dots,\beta_m)\in [p_2]^m$. We say that $\pi$ is
{\it order-isomorphic\/} to $\beta$ if for all $1\leq i<j\leq m$ one has
$\pi_i<\pi_j$ if and only if $\beta_i<\beta_j$. For two permutations
$\pi\in S_n$ and $\tau\in S_k$, an {\it occurrence\/} of $\tau$ in $\pi$
is a subsequence $1\leq i_1<i_2<\dots<i_k\leq n$ such that $(\pi_{i_1},
\dots,\pi_{i_k})$ is order-isomorphic to $\tau$; in such a context $\tau$ is
usually called the {\it pattern \/} (classical pattern). 
We say that $\pi$ {\it avoids\/} $\tau$,
or is $\tau$-{\it avoiding\/}, if there is no occurrence of $\tau$ in $\pi$. 
More generally, we say $\pi$ {\em containing} $\tau$ exactly $r$ times, if 
there exists $r$ different occurrences of $\tau$ in $\pi$. \\
The set of all $\tau$-avoiding permutations of all possible sizes including 
the empty permutation is denoted $\mathcal{S}(\tau)$. Pattern avoidance proved
to be a useful language in a variety of seemingly unrelated problems, from
stack sorting \cite{Kn} to singularities of Schubert varieties \cite{LS}.
A complete study of pattern avoidance for the case $\tau\in S_3$ is
carried out in \cite{SS}. \\

On the other hand, in \cite{BS} introduced generalised permutation 
patterns that allow the requirement that two adjacent letters in a pattern
must be adjacent in the permutation. The idea of \cite{BS} 
introducing these patterns was study of Mahonian statistics. 

We write a classical pattern with dashes between any two adjacent letters of the pattern, 
say $1324$, as $1\mn3\mn2\mn4$, and if we write, say $24\mn1\mn3$, then we mean that 
if this pattern occurs in permutation $\pi$, then the letters in the 
permutation $\pi$ that correspond to $2$ and $4$ are adjacent.
For example, the permutation $\pi=35421$ has only two occurrences of the pattern $23\mn1$, 
namely the subsequnces $352$ and $351$, whereas $\pi$ has four occurrences of the 
pattern $2\mn3\mn1$, namely the subsequnces $352$, $351$, $342$ and $341$. 

In \cite{C} presented a complete solution for the number of permutations avoiding any 
pattern of length three with exactly one adjacent pair of letters. In \cite{CM} 
presented a complete solution for the number of permutations avoiding any two patterns of 
length three with exactly one adjacent pair of letters. In \cite{Ki} 
presented almost results avoiding two or more $3$-patterns without internal dashes, that is, 
where the pattern corresponds to a contiguous subword in a permutation. Besides, \cite{EN} 
presented the following generating functions regarding the distribution of the number of 
occurrences of any generalized pattern of length $3$:
$$\begin{array}{ll}
\sum\limits_{\pi\in\mathcal{S}} y^{(123)\pi} \frac{x^{|\pi|}} {|\pi|!}&=
	\frac{2f(y)e^{\frac{1}{2}(f(y)-y+1)x}} {f(y)+y+1+(f(y)-y-1)e^{f(y)x}},\\
\sum\limits_{\pi\in\mathcal{S}} y^{(213)\pi} \frac{x^{|\pi|}}{|\pi|!}&=
	\frac{1}{1-\int\limits_0^x e^{(y-1)t^2/2}dt},
\end{array}$$
where $(\tau)\pi$ is the number of occurrences of $\tau$ in $\pi$, $f(y)=\sqrt{(y-1)(y+3)}$.\\

The purpose of this paper is to point out an analog of a theorem \cite{RWZ}, and some 
interesting consequences of this analog. Generalization of this theorem have already 
been given in \cite{Kr,MV1,JR}. 
In the present note we study the generating function for the number 
$1$-$3$-$2$-avoiding permutations in $S_n$ that contain a prescribed number of generalised 
pattern $\tau$. The study of the obtained continued fraction allows us to recover and to present 
analog of the results of \cite{CW,Kr,MV1,JR} that relates the  number of 
$1$-$3$-$2$-avoiding permutations that contain no $12$-$3$-$\dots$-$k$ 
(or $21$-$3$-$\dots$-$k$)
patterns to Chebyshev polynomials of the second kind. \\

Let $f_{\tau;r}(n)$ stand for the number of $1$-$3$-$2$-avoiding 
permutations in $S_n$ that contain exactly $r$ occurrences of $\tau$. 
We denote by $F_{\tau;r}(x)$ and $F_\tau(x,y)$ the generating function 
of the sequence $\{f_{\tau;r}(n)\}_{n\geq 0}$ and $\{f_{\tau;r}(n)\}_{n,r\geq 0}$
respectively, that is,
	$$F_{\tau;r}(x)=\sum_{n\geq0} f_{\tau;r}(n)x^n,\quad F_{\tau}(x,y)=\sum_{r\geq 0} F_{\tau;r}(x)y^r.$$

The paper is organized as follows. The cases $\tau=12\mn3\mn\dots\mn k$, 
$\tau=21\mn3\mn\dots\mn k$, $\tau=123\dots k$, and $\tau=k\dots 321$ are 
treat in {\bf section 2}. 
In {\bf section $3$}, we treat the cases $\tau=123$, $213$, $231$, $312$, and $321$, 
that is $\tau$ is $3$-letters generalized pattern without dashes. In {\bf section $4$}, 
we treat the cases when $\tau$ is $3$-lettesr generalized pattern with one dash. Finally,  
in {\bf section $5$}, we treat some of examples of restricted more than one generalized 
pattern of $3$-letters. 
%=====================================================================
\section{FOUR GENERAL CASES}
In this section, we study the following four cases $\tau=12\mn3\mn\dots\mn k$, 
$\tau=21\mn3\mn\dots\mn k$, $\tau=12\dots k$, 
and $\tau=k\dots 21$, by the following three subsections.

%*************
\section*{2.1 Pattern $12\mn3\mn\dots\mn k$}
Our first result is a natural analog of the main theorems of \cite{RWZ,MV1,Kr}.

\begin{theorem} 
\label{th21}
The generating function $F_{12\mn3\mn\dots\mn k}(x,y)$ for $k\geq2$ is given by the continued fraction
$$
\frac1{1-x+xy^{d_1}-\dfrac{xy^{d_1}}{1-x+xy^{d_2}-\dfrac{xy^{d_2}}
{1-x+xy^{d_3}-\dfrac{xy^{d_3}}{\dots}}}},
$$
where $d_i=\binom{i-1}{k-2}${\rm,} and $\binom ab$ is assumed $0$ whenever
$a<b$ or $b<0$.
\end{theorem}
\begin{proof} Following \cite{MV1} we define $\eta_j(\pi)$, $j\geq3$, as the
number of occurrences of $12\mn3\mn\dots\mn j$ in $\pi$. Define $\eta_2(\pi)$ for any $\pi$, 
as the number of occurrences of $12$ in $\pi$, $\eta_1(\pi)$ as the number letters 
of $\pi$, and $\eta_0(\pi)=1$ for any $\pi$, which means that the empty pattern 
occurs exactly once in each permutation. The {\it weight\/} of a permutation 
$\pi$ is a monomial in $k$ independent variables $q_1,\dots,q_k$ defined by
	$$w_k(\pi)=\prod_{j=1}^k q_j^{\eta_j(\pi)}.$$
The {\it total weight\/} is a polynomial
	$$W_k(q_1,\dots,q_k)=\sum_{\pi\in\sss}w_k(\pi).$$
The following proposition is implied immediately by the definitions.

\begin{proposition} 
\label{pro21}
$F_{12\mn3\mn\dots\mn k}(x,y)=W_k(x,1,\dots,1,y)$ for $k\geq2$.
\end{proposition}

We now find a recurrence relation for the numbers $\eta_j(\pi)$.
Let $\pi\in S_n$, so that $\pi=(\pi',n,\pi'')$.

\begin{proposition} 
\label{pro22}
For any nonempty $\pi\in\sss$
	$$\eta_j(\pi)=\eta_j(\pi')+\eta_j(\pi'')+\eta_{j-1}(\pi'),$$
where $j\neq 2$. Besides, if $\pi'$ nonempty then 
	$$\eta_2(\pi)=\eta_2(\pi')+\eta_2(\pi'')+1,$$
otherwise 
	$$\eta_2(\pi)=\eta_2(\pi'').$$
\end{proposition}
\begin{proof} Let $l=\pi^{-1}(n)$.
Since $\pi$ avoids $1\mn3\mn2$, each number in $\pi'$ is greater than
any of the numbers in $\pi''$. Therefore, $\pi'$ is a $1-3-2$-avoiding
permutation of the numbers $\{n-l+1,n-l+2,\dots,n-1\}$, while $\pi''$
is a $1\mn3\mn2$-avoiding permutation of the numbers $\{1,2,\dots,n-l\}$. On the
other hand, if $\pi'$ is an arbitrary $1\mn3\mn2$-avoiding
permutation of the numbers $\{n-l+1,n-l+2,\dots,n-1\}$ and $\pi''$ is an
arbitrary $1\mn3\mn2$-avoiding permutation of the numbers $\{1,2,\dots,n-l\}$, then
$\pi=(\pi',n,\pi'')$ is $1-3-2$-avoiding. Finally, if $(i_1,\dots,i_j)$ is
an occurrence of $12\mn3\mn\dots\mn j$ in $\pi$ then either $i_j<l$, and so it is
also an occurrence  of $12\mn3\mn\dots\mn j$ in $\pi'$, or $i_1>l$, and so it is
also an occurrence  of $12\mn\dots\mn j$ in $\pi''$, or $i_j=l$, and so
$(i_1,\dots,i_{j-1})$ is an occurrence  of $12\mn3\mn\dots\mn(j-1)$ in $\pi'$, where 
$j\neq 2$. For $j=2$ the proposition is trivial. The result follows.
\end{proof}

Now we are able to find the recurrence relation for the total weight $W$.
Indeed, by Proposition \ref{pro22},
$$
\begin{array}{lll}
W_k(q_1,\dots,q_k)&=1+\sum\limits_{\varnothing\ne\pi\in\sss}&\prod\limits_{j=1}^k   
        q_j^{\eta_j(\pi)}\\
&=1+\sum\limits_{\emptyset\neq\pi'\in\sss} &\sum\limits_{\pi''\in\sss} \prod\limits_{j=1}^k
q_j^{\eta_j(\pi'')}\cdot q_1^{\eta_1(\pi')+1}q_2\cdot\\
&	& \prod\limits_{j=2}^{k-1}(q_jq_{j+1})^{\eta_j(\pi')}\cdot q_k^{\eta_k(\pi')} +\sum\limits_{\pi''\in\sss} q_1\prod\limits_{j=1}^k q_j^{\eta_j(\pi'')}.
\end{array}
$$
Hence
$$\begin{array}{ll}
W_k(q_1,\dots,q_k)&=1+q_1W_k(q_1,\dots,q_k)+\\
		    &+q_1q_2W_k(q_1,\dots,q_k)(W_k(q_1, q_2q_3,\dots,q_{k-1}q_k,q_k)-1).
\end{array} \eqno(1)$$

Following \cite{MV1}, for any $d\geq0$ and $2\leq m\leq k$ define 
	$$\qq^{d,m}=\prod_{j=2}^k q_j^{{{d}\choose{j-m}}};$$
recall that ${a\choose b}=0$ if $a<b$ or $b<0$. The following proposition is
implied immediately by the well-known properties of binomial coefficients.

\begin{proposition}
\label{pro23} 
For any $d\geq0$ and $2\leq m\leq k$
	$$\qq^{d,m}\qq^{d,m+1}=\qq^{d+1,m}.$$
\end{proposition}

Observe now that $W_k(q_1,\dots,q_k)=W_k(q_1,\qq^{0,2},\dots,\qq^{0,k})$ and that
by (1) and Proposition \ref{pro23}

$$\begin{array}{ll}
W_k(q_1,\qq^{d,2},\dots,\qq^{d,k})&=1+q_1W_k(q_1,\qq^{d,2},\dots,\qq^{d,k})+\\
		    &+q_1q^{d,2}W_k(q_1,q^{d,2},\dots,\qq^{d,k})(W_k(q_1,\qq^{d+1,2},\dots,\qq^{d+1,k})-1).
\end{array}$$
therefore
	$$W_k(q_1,\dots,q_k)=\frac1{1-q_1+q_1\qq^{0,2}-\dfrac{q_1\qq^{0,2}}{1-q_1+q_1\qq^{1,2}-\dfrac{q_1\qq^{1,2}}
	{1-q_1+q_1\qq^{2,2}-\dfrac{q_1\qq^{2,2}}{\dots}}}}.$$
To obtain the continued fraction representation for $F(x,y;k)$ it is enough
to use Proposition \ref{pro21} and to observe that
	$$q_1\qq^{d,2}\bigg|_{q_1=x,q_2=\dots=q_{k-1}=1,q_k=y}=xy^{\binom{d}{k-2}}.$$
\end{proof} 

\begin{corollary} 
\label{c25}
	$$F(x,y;2)=\frac{1-x+xy-\sqrt{(1-x)^2-2x(1+x)y+x^2y^2}}{2xy},$$
in other words, for any $r\geq 1$
	$$f_{12}(n)=\frac{r+1}{n(n-r)}{{n}\choose{r+1}}^2.$$
\end{corollary}
\begin{proof}
For $k=2$, Theorem \ref{th21} yields
$$
F_{12}(x,y)=\frac1{1-x+xy-\dfrac{xy}{1-x+xy-\dfrac{xy}
{1-x+xy-\dfrac{xy}{\dots}}}}. 
$$
which means that 
	$$F_{12}(x,y)=\frac{1}{1-x+xy-xyF_{12}(x,y)}.$$
So the rest it is easy to see.
\end{proof}

Now, we find explicit expressions for $F_{12\mn3\mn\dots\mn k;r}(x)$ 
in the case $0\leq r\leq k-2$.

Following \cite{MV1}, consider a recurrence relation
	$$T_j=\frac1{1-xT_{j-1}},\quad j\gs1. \eqno(2)$$
The solution of (2) with the initial condition $T_0=0$ is denoted by
$R_j(x)$, and the solution of (2) with the initial condition
$$
T_0=G_{12\mn3\mn\dots\mn k}(x,y)=\frac 1{1-x+xy^{\binom {k-2}0}-\dfrac{xy^{\binom {k-2}0}}{1-x+xy^{\binom{k-1}1}-\dfrac{xy^{\binom{k-1}1}}
{1-x+xy^{\binom{k}2}-\dfrac{xy^{\binom{k}2}}{\dots}}}}
$$
is denoted by $S_j(x,y;k)$, or just $S_j$ when the value of $k$ is clear from
the context. Our interest in (2) is stipulated by the
following relation, which is an easy consequence of Theorem \ref{th21}:
	$$F_{12\mn3\mn\dots\mn k}(x,y)=S_{k-2}(x,y;k). \eqno(3)$$
Following to \cite[eq. 4]{MV1}, for all $j\geq 1$
	$$R_j(x)=\frac{U_{j-1}\left(\frac1{2\sqrt{x}}\right)}
	{\sqrt{x}U_{j}\left(\frac1{2\sqrt{x}}\right)}, 
	\eqno(4)$$
where $U_j(\cos\theta)=\sin(j+1)\theta/\sin\theta$ be the Chebyshev polynomials
of the second kind. Next, we find an explicit expression for $S_j$ in terms of $G$
and $R_j$.

\begin{lemma}
\label{lem32} 
For any $j\geq 2$ and any $k\geq 2$
$$
S_j(x,y;k)=R_j(x)\frac{1-xR_{j-1}(x)G_{12\mn3\mn\dots\mn k}(x,y)}{1-xR_{j}(x)G_{12\mn3\mn\dots\mn k}(x,y)}. \eqno(5)
$$
\end{lemma}

\begin{proof} Indeed, from (2) and $S_0=G$ we get
$S_1=1/(1-xG)$. On the other hand, $R_0=0$, $R_1=1$, so (5)
holds for $j=1$. Now let $j>1$, then by induction
$$
S_j=\frac1{1-xS_{j-1}}=\frac1{1-xR_{j-1}}\cdot
\frac{1-xR_{j-1}G}{1-\dfrac{x(1-xR_{j-2})R_{j-1}G}
{1-xR_{j-1}}}.
$$
Relation (2) for $R_j$ and $R_{j-1}$ yields
$(1-xR_{j-2})R_{j-1}=(1-xR_{j-1})R_j=1$, which together with the above
formula gives (5).
\end{proof}

As a corollary from Lemma \ref{lem32} and (3) we get the following expression for
the generating function $F_{12\mn3\mn\dots\mn k}(x,y)$.

\begin{corollary} For any $k\geq 3$
$$
F_{12\mn3\mn\dots\mn k}(x,y)=R_{k}(x)+\big(R_{k-2}(x)-R_{k-3}(x)\big)
\sum_{m\geq 1}\big(xR_{k-2}(x)G_{12\mn3\mn\dots\mn k}(x,y)\big)^m.
$$
\end{corollary}

Now we are ready to express the generating functions $F_{12\mn3\mn\dots\mn k;r}(x)$,
$0\ls r\ls k-2$, via Chebyshev polynomials.

\begin{theorem} 
\label{th31} For any $k\geq 3${\rm,} $F_{12\mn3\mn\dots\mn k;r}(x)$ 
is a rational function given by
$$
\begin{array}{ll}
F_{12\mn3\mn\dots\mn k;r}(x)&=\frac{x^{r-1}U_{k-2}^{r-1}\left(\frac1{2\sqrt{x}}\right)}
{(1-x)^rU_{k}^{r+1}\left(\frac1{2\sqrt{x}}\right)},\quad 1\leq r\leq k-2,\\
F_{12\mn3\mn\dots\mn k;0}(x)&=\frac{U_{k-1}\left(\frac1{2\sqrt{x}}\right)}
{\sqrt{x}U_{k}\left(\frac1{2\sqrt{x}}\right)},
\end{array}
$$
where $U_j$ is the $j$th Chebyshev polynomial of the second kind.
\end{theorem}
\begin{proof} Observe that $G_{12\mn3\mn\dots\mn k}(x,y)=\frac{1}{1-x}\cdot\frac{1}{1-\frac{x^2}{(1-x)^2}y}+y^{k-1}P(x,y)$,
so from Corollary we get
$$\begin{array}{ll}
F_{12\mn3\mn\dots\mn k}(x,y)=R_k(x)&+\big(R_{k-2}(x)-R_{k-3}(x)\big)
\sum\limits_{m=1}^k\big(\frac{x}{1-x}R_k(x)\big)^m\sum\limits_{n=1}^{k-2} {\binom {m-1+n}n} \frac{x^{2n}}{(1-x)^{2n}} y^n\\
&+y^{k-1}P'(x,y),
\end{array}$$
where $P(x,y)$ and $P'(x,y)$ are formal power series.
To complete the proof, it suffices to use (4) together with the identity
$U_{n-1}^2(z)-U_n(z)U_{n-2}(z)=1$,
which follows easily from the trigonometric identity
$\sin^2n\theta-\sin^2\theta=\sin(n+1)\theta\sin(n-1)\theta$.
\end{proof}

%******************
\section*{2.2 Pattern $21\mn3\mn\dots\mn k$}
Our second result is a natural analog of the main theorems of \cite{RWZ,MV1,Kr}.

\begin{theorem}
\label{th21b}
For any $k\geq 2$,
$$
F_{21\mn3\mn\dots\mn k}(x,y)=
1-\frac{x}
{xy^{d_1}-\dfrac{1}
{1-\dfrac{x}
{xy^{d_2}-\dfrac{1}
{1-\dfrac{x}
{xy^{d_3}-\dfrac{1}{\ddots}}}}}},$$
where $d_i=\binom{i-1}{k-2}${\rm,} and $\binom ab$ is assumed $0$ whenever
$a<b$ or $b<0$.
\end{theorem}
\begin{proof} Following \cite{MV1} we define $\nu_j(\pi)$, $j\geq3$, as the
number of occurrences of $21\mn3\mn\dots\mn j$ in $\pi$. Define $\nu_2(\pi)$ for any $\pi$, 
as the number of occurrences of $21$ in $\pi$, $\nu_1(\pi)$ as the number letters 
of $\nu$, and $\nu_0(\pi)=1$ for any $\pi$, which means that the empty pattern 
occurs exactly once in each permutation. The {\it weight\/} of a permutation 
$\pi$ is a monomial in $k$ independent variables $q_1,\dots,q_k$ defined by
	$$v_k(\pi)=\prod_{j=1}^k q_j^{\nu_j(\pi)}.$$
The {\it total weight\/} is a polynomial
	$$V_k(q_1,\dots,q_k)=\sum_{\pi\in\sss}v_k(\pi).$$
The following proposition is implied immediately by the definitions.

\begin{proposition} 
\label{pro21b}
$F_{21\mn3\mn\dots\mn k}(x,y)=V_k(x,1,\dots,1,y)$ for $k\geq2$.
\end{proposition}

We now find a recurrence relation for the numbers $\nu_j(\pi)$.
Let $\pi\in S_n$, so that $\pi=(\pi',n,\pi'')$.

\begin{proposition} 
\label{pro22b}
For any nonempty $\pi\in\sss$
	$$\nu_j(\pi)=\nu_j(\pi')+\nu_j(\pi'')+\nu_{j-1}(\pi'),$$
where $j\neq 2$. Besides, if $\pi''$ nonempty then 
	$$\nu_2(\pi)=\nu_2(\pi')+\nu_2(\pi'')+1,$$
otherwise 
	$$\nu_2(\pi)=\nu_2(\pi'').$$
\end{proposition}
\begin{proof} Similarly to Proposition \ref{pro22} we get  
$\pi$ avoids $1\mn3\mn2$ if and only if $\pi'$ is a $1\mn3\mn2$-avoiding
permutation of the numbers $\{n-l+1,n-l+2,\dots,n-1\}$, while $\pi''$
is a $1-3-2$-avoiding permutation of the numbers $\{1,2,\dots,n-l\}$. 
Finally, if $(i_1,\dots,i_j)$ is
an occurrence of $21\mn3\mn\dots\mn j$ in $\pi$ then either $i_j<l$, and so it is
also an occurrence  of $21\mn3\mn\dots\mn j$ in $\pi'$, or $i_1>l$, and so it is
also an occurrence  of $21\mn3\mn\dots\mn j$ in $\pi''$, or $i_j=l$, and so
$(i_1,\dots,i_{j-1})$ is an occurrence  of $21\mn3\mn\dots\mn (j-1)$ in $\pi'$, where 
$j\neq 2$. For $j=2$ the proposition is trivial. The result follows.
\end{proof}

Now we are able to find the recurrence relation for the total weight $V$.
Indeed, by Proposition \ref{pro22b},
$$
\begin{array}{lll}
V_k(q_1,\dots,q_k)&=1+\sum\limits_{\varnothing\ne\pi\in\sss}&\prod\limits_{j=1}^k   
        q_j^{\nu_j(\pi)}\\
&=1+\sum\limits_{\emptyset\neq\pi''\in\sss} &\sum\limits_{\pi'\in\sss} \prod\limits_{j=1}^k
q_j^{\nu_j(\pi'')}\cdot  q_1^{\nu_1(\pi')+1}q_2\cdot\prod\limits_{j=2}^{k-1}(q_jq_{j+1})^{\nu_j(\pi')}\cdot q_k^{\nu_k(\pi')}+\\
& &+\sum\limits_{\pi'\in\sss} q_1q_1^{\nu(\pi')}q_k^{\nu_k(\pi')}\prod\limits_{j=2}^{k-1} (q_jq_{j+1})^{\nu_j(\pi')}.
\end{array}
$$
Hence
$$\begin{array}{ll}
V_k(q_1,\dots,q_k)&=1+q_1V_k(q_1,q_2q_3,\dots,q_{k-1}q_k,q_k)+\\
		    &+q_1q_2V_k(q_1, q_2q_3,\dots,q_{k-1}q_k,q_k)(V_k(q_1,q_2,\dots,q_k)-1). 
\end{array}\eqno(6)$$

Observe now that $V_k(q_1,\dots,q_k)=V_k(q_1,\qq^{0,2},\dots,\qq^{0,k})$ and that
by (6) and Proposition \ref{pro23}

$$\begin{array}{ll}
V_k(q_1,\qq^{d,2},\dots,\qq^{d,k})&=1+q_1V_k(q_1,\qq^{d+1,2},\dots,\qq^{d+1,k})+\\
		    &+q_1q^{d,2}V_k(q_1,q^{d+1,2},\dots,\qq^{d+1,k})(V_k(q_1,\qq^{d,2},\dots,\qq^{d,k})-1).
\end{array}$$

To obtain the continued fraction representation for $F_{21\mn3\mn\dots\mn k}(x,y)$ it is enough
to use Proposition \ref{pro21b} and to observe that
	$$q_1\qq^{d,2}\bigg|_{q_1=x,q_2=\dots=q_{k-1}=1,q_k=y}=xy^{\binom{d}{k-2}}.$$
\end{proof}

\begin{corollary} 
\label{c25b}
	$$F_{21}(x,y)=\frac{1-x+xy-\sqrt{(1-x)^2-2x(1+x)y+x^2y^2}}{2xy},$$
in other words, for any $r\geq 1$
	$$f_{21;r}(n)=\frac{r+1}{n(n-r)}{{n}\choose{r+1}}^2.$$
\end{corollary}
\begin{proof}
For $k=2$, $q_1=x$ and $q_2=y$, Proposition \ref{pro21b} and (6) yields
$F_{21}(x,y)=1+xF_{21}(x,y)+xyF_{21}(x,y)(F_{21}(x,y)-1)$, which means $F_{21}(x,y)=F_{12}(x,y)$. 
The rest of the proof obtained by Corollary \ref{c25}.
\end{proof}

Now, we are ready to find an 
explicit expressions for $F_{21\mn3\mn\dots\mn k;r}(x)$ in the case $0\leq r\leq k-2$.\\

Consider a recurrence relation
	$$T'_j=1-\frac{x}{x-\dfrac{1}{T'_{j-1}(x)}},\quad j\geq1. \eqno(7)$$
The solution of (7) with the initial condition $T'_0=0$ is given by
$R_j(x)$ (Lemma \ref{lem32b}), and the solution of (7) with the initial condition
$$
T'_0=G_{21\mn3\mn\dots\mn k}(x,y)=
\frac{1}
{xy^{\binom{k-2}{k-2}}-\dfrac{1}
{1-\dfrac{x}
{xy^{\binom{k-1}{k-2}}-\dfrac{1}
{1-\dfrac{x}
{xy^{\binom{k}{k-2}}-\dfrac{1}{\ddots}}}}}},
$$
is denoted by $S'_j(x,y;k)$, or just $S'_j$ when the value of $k$ is clear from
the context. Our interest in (7) is stipulated by the
following relation, which is an easy consequence of Theorem \ref{th21b}:
		$$F_{21\mn3\mn\dots\mn k}(x,y)=S'_k(x,y;k). \eqno(8)$$

First of all, we find an explicit formula for the functions $T'_j(x)$ in (7).

\begin{lemma} 
\label{lem32b} For any $j\geq 1$, 
	$$T'_j(x)=R_j(x). \eqno(9)$$
\end{lemma}
\begin{proof} Indeed, it follows immediately from (7) that $T'_0(x)=0$ 
and $T'_1(x)=1$. Let us induction, we assume $T'_{j-1}(x)=R_{j-1}(x)$, and prove 
that $T'_j(x)=R_j(x)$. By use (7)
	$$T'_j(x)=1-\frac{x}{x-\dfrac{1}{R_{j-1}(x)}}.$$
On the other hand, following to \cite{MV1}, $R_j(x)=\frac{1}{1-xR_{j-1}(x)}$ 
which means that $R_j(x)=1+xR_{j-1}(x)R_j(x)$, hence $T'_j(x)=R_j(x)$.
\end{proof}

Next, we find an explicit expression for $S'_j$ in terms of $G$ and $R_j$.

\begin{lemma} 
\label{lem33b} For any $j\geq 2$ and any $k\geq 2$
$$
S'_j(x,y;k)=R_j(x)\frac{1-xR_{j-1}(x)G_{21\mn3\mn\dots\mn k}(x,y;k)}{1-xR_{j}(x)G_{21\mn3\mn\dots\mn k}(x,y)}. \eqno(10)
$$
\end{lemma}

As a corollary from Lemma \ref{lem33b} and (6) we get the following expression for
the generating function $F_{21\mn3\mn\dots\mn k}(x,y)$.

\begin{corollary} For any $k\geq 3$
\label{cc1}
$$F_{21\mn3\mn\dots\mn k}(x,y)=R_{k}(x)+\big(R_{k-2}(x)-R_{k-3}(x)\big)\sum_{m\geq 1}\big(xR_{k-2}(x)G_{21\mn3\mn\dots\mn k}(x,y)\big)^m.$$
\end{corollary}

Now we are ready to express the generating functions $F_{21\mn3\mn\dots\mn k;r}(x)$,
$0\ls r\ls k-2$, via Chebyshev polynomials.

\begin{theorem} 
\label{th31b} For any $k\geq 3${\rm,} $F_{21\mn3\mn\dots\mn k;r}(x)$ is a rational function
given by
$$
\begin{array}{ll}
F_{21\mn3\mn\dots\mn k;r}(x)&=\frac{x^{\frac{r-1}2}U_{k-2}^{r-1}\left(\frac1{2\sqrt{x}}\right)}
{U_{k}^{r+1}\left(\frac1{2\sqrt{x}}\right)},\quad 1\leq r\leq k-2,\\
F_{21\mn3\mn\dots\mn k;0}(x)&=\frac{U_{k-1}\left(\frac1{2\sqrt{x}}\right)}
{\sqrt{x}U_{k}\left(\frac1{2\sqrt{x}}\right)},
\end{array}
$$
where $U_j$ is the $j$th Chebyshev polynomial of the second kind.
\end{theorem}
\begin{proof} Observe that $G_{21\mn3\mn\dots\mn k}(x,y)=1+\frac{x}{1-x-xy}+y^{k-1}P(x,y)$,
so from Corollary \ref{cc1} we get
$$
F_{21\mn3\mn\dots\mn k}(x,y)=R_k(x)+\big(R_{k-2}(x)-R_{k-3}(x)\big)
\sum\limits_{m=1}^k\big(xR_{k-2}(x)\big( 1+ \frac{x}{1-x-xy} \big) \big)^m+y^{k-1}P'(x,y),
$$
where $P(x,y)$ and $P'(x,y)$ are formal power series.
To complete the proof, it suffices to use (9) together with the identity
$U_{n-1}^2(z)-U_n(z)U_{n-2}(z)=1$.
\end{proof}

\begin{remark}
Theorem \ref{th31b} and \cite{MV1} yields the number of $1\mn3\mn2$-avoiding 
permutations in $S_n$ such that containing $r$ times the pattern $21\mn3\mn\dots\mn k$ is the same 
number of $1\mn3\mn2$-avoiding permutations in $S_n$ such that containing $r$ times the pattern 
$1\mn2\mn3\dots\mn k$, for all $r=0,1,2,\dots,k-2$. However, the question if there exist 
a natural bijection between the set of $1\mn3\mn2$-avoiding permutations in $S_n$ such containing $r$ 
times the generalized pattern $21\mn3\mn\dots\mn k$, and the set of $1\mn3\mn2$-avoiding 
permutations in $S_n$ such containing $r$ times the classically pattern $1\mn 2\mn3\mn\dots\mn k$.
\end{remark}
%*********************************************************
\section*{2.3 Patterns: $\tau=12\dots k$ and $\tau=k\dots 21$}
Let $\pi\in S_n$; we say $\pi$ has {\it $d$-increasing canonical decomposition\/} 
if $\pi$ has the following form
	$$\pi=(\pi^1,\pi^2,\dots,\pi^d,a_d,\dots,a_2,a_1,n,\pi^{d+1}),$$
where all the entries of $\pi^i$ are greater than all the entries of $\pi^{i+1}$, and 
$a_d<a_{d-1}<\dots<a_1<n$. We say $\pi$ has {\it $d$-decreasing canonical decomposition\/} 
if $\pi$ has the following form
	$$\pi=(\pi^1,n,a_1,\dots,a_d,\pi^{d+1},\pi^d,\dots,\pi^d),$$
where all the entries of $\pi^i$ are greater than all the entries of $\pi^{i+1}$, and 
$a_d<a_{d-1}<\dots<a_1<n$.
The following proposition it the base of all the other
results in this Section.

\begin{proposition}
\label{progen}
Let $\pi\in S_n(1\mn3\mn2)$. Then there exists unique $d\geq 0$ and $e\geq 0$ such that 
$\pi$ has a $d$-increasing canonical decomposition, and has $e$-decreasing canonical decomposition.
\end{proposition}
\begin{proof}
Let $\pi\in S_n(1\mn3\mn2)$, and let $a_d,a_{d-1},\dots,a_1,n$ a maximal increasing 
subsequence of $\pi$ such that $\pi=(\pi',a_d,\dots,a_1,n,\pi'')$. 
Since $\pi$ avoids $1\mn3\mn2$ there exists 
$d$ subsequnces $\pi^j$ such that $\pi=(\pi^1,\dots,\pi^d,a_d,\dots,a_1,n,\pi'')$, and 
all the entries of $\pi^i$ are greater than all the entries of $\pi^{i+1}$, and all the 
entries of $\pi^d$ greater than all entries of $\pi''$. Hence, 
$\pi$ has $d$-increasing canonical decomposition. Similarly, there exist $e$ unique 
such that $\pi$ is $e$-decreasing canonical decomposition.
\end{proof}

Let us define $I_\tau(x,y;d)$ (respectively; $J_\tau(x,y;e)$) be the generating function for all 
$d$-increasing (respectively; $e$-decreasing) canonical decomposition of permutations 
in $S_n(1\mn3\mn2)$ with exactly $r$ occurrences of $\tau$. 
The following proposition is implied immediately by the definitions.

\begin{proposition}
\label{progen1}
$$F_\tau(x,y)=1+\sum_{d\geq 0} I_\tau(x,y;d)=1+\sum_{e\geq 0} J_\tau(x,y;e).$$
\end{proposition}
\begin{proof}
Immetaitley, by definitions of the generating functions and 
Proposition \ref{progen} ($1$ for the empty permutation).
\end{proof}

Now, we present an examples for Proposition \ref{progen} 
and Proposition \ref{progen1}.

\subsubsection*{First example: $F_{12\dots k}(x)$ and $F_{k\dots 21}(x)$}

\begin{theorem}
\label{th21c}
$F_{k\dots 21}(x,y)=F_{12\dots k}(x,y),$ such that
	$$F_{12\dots k}(x,y)=\sum_{n=0}^{k-2} x^nF^n_{12\dots k}(x,y)+\frac{x^{k-1}F^{k-1}_{12\dots k}(x,y)}{1-xyF_{12\dots k}(x,y)}.$$
\end{theorem}
\begin{proof}
By Proposition \ref{progen} and definitions it is easy to obtain for all $d\geq 0$
	$$I_{12\dots k}(x,y;d)=x^{d+1}y^{s_d}F_{12\dots k}^{d+1}(x,y),$$
where $s_d=d+1-k$ for $d\geq k-1$, and otherwise $s_d=0$. So by Proposition \ref{progen1} 
the theorem holds.

Similarly, we obtain the same result for $F_{k\dots 21}(x,y)$. 
\end{proof}

As a remark, by the above theorem, it is easy to obtain the same 
Corollaries \ref{c25} and \ref{c25b} results.

\subsubsection*{Second example}
Here we calculate $F_{1\mn2\mn\dots\mn (l-1)\mn l(l+1)\dots k}(x,y)$ where $l\leq k-1$. 
\begin{theorem}
\label{thgen1}
Let $1\leq l\leq k-1$. Then 
$F_{1\mn2\mn\dots\mn (l-1)\mn l(l+1)\dots k}(x,y)=U_l(x,1,\dots,1,y)$ where
$$U_l(q_1,\dots,q_l)=1+\sum_{d\geq 0} \left( (q_l^{\binom{d+1+l-k}{l}}\prod_{j=1}^{l-1} q_j^{\binom{d+1}{j}} \prod_{j=0}^d U_l(p_{1;j},\dots,p_{l;j}) \right),$$
and for $i=1,2,\dots,l$,  
$p_{i;j}=\prod_{m=1}^{l-1} q_j^{\binom{j}{m-i}}$, $p_{l,j}=q_l$ for all $0\leq j\leq k-l$, 
and $p_{i;j}=\prod_{m=1}^l p_{i;k-l}^{\binom{j-k+l}{l-i}}$ for all $j\geq k-l+1$.
\end{theorem}
\begin{proof}
Following \cite{MV1} we define $\gamma_j(\pi)$, $j\leq l-1$, as the
number of occurrences of $1\mn2\mn\dots\mn j$ in $\pi$. Define $\gamma_l(\pi)$ for any $\pi$, 
as the number of occurrences of $1\mn2\mn\dots\mn (l-1)\mn l(l+1)\dots k$ in $\pi$, and 
$\gamma_0(\pi)=1$ for any $\pi$, which means that the empty pattern 
occurs exactly once in each permutation. The {\it weight\/} of a permutation 
$\pi$ is a monomial in $l$ independent variables $q_1,\dots,q_l$ defined by
	$$u_l(\pi)=\prod_{j=1}^l q_j^{\gamma_j(\pi)}.$$
The {\it total weight\/} is a polynomial
	$$U_l(q_1,\dots,q_l)=\sum_{\pi\in\sss}u_l(\pi).$$
The following proposition is implied immediately by the definitions, 
and Proposition \ref{progen}.

\begin{proposition} 
\label{pro211}
$F_{1\mn2\mn\dots\mn (l-1)\mn l(l+1)\dots k}(x,y)=U_k(x,1,\dots,1,y)$ for $k>l\geq1$, and 
$U_l(q_1,\dots,q_l)=1+\sum_{d\geq 0}\sum_{\pi\in A_d} u_l(\pi)$, 
where $A_d$ is the set of all $d$-increasing canonical decomposition 
permutations in $\sss$.
\end{proposition}

Let us denote $U_{l;d}(q_1,\dots,q_l)=\sum_{\pi\in A_d} u_l(\pi)$. 

\begin{proposition}
\label{pro212} For any $d\geq 0$, 
$$U_{l;d}(q_1,\dots q_l)=q_l^{\binom{d+1+l-k}{l}}\prod_{j=1}^{l-1} q_j^{\binom{d+1}{j}} \prod_{j=0}^d U_l(p_{1;j},\dots,p_{l;j}).$$
\end{proposition}
\begin{proof}
Let $\pi$ is $d$-increasing canonical decomposition; that is,  
	$$\pi=(\pi^1,\pi^2,\dots,\pi^d,a_d,\dots,a_2,a_1,n,\pi^{d+1}),$$
where the numbers $a_{d}<a_{d-1}<\dots<a_1<n$ are appear as 
consecutive numbers in $\pi$, the all entries of $\pi^j$ are greater than 
all the entries of $\pi^{j+1}$, and the all entries of $\pi^d$ greater than 
$a_d$. So, by calculate $u_l(\pi)$ and sum over all $\pi\in A_d$ we 
have that 
  $$U_{l;d}(q_1,\dots,q_d)=q_l^{\binom{d+1+l-k}{l}} 
     	\cdot \prod_{j=1}^{l-1} q_j^{\binom{d+1}{j}} 
	\cdot \prod_{j=0}^d U_l(p_{1;j},\dots,p_{l;j}).$$
\end{proof}

Therefore, Theorem \ref{thgen1} holds, by use the above results; 
Proposition \ref{pro211} and Proposition \ref{pro212}.
\end{proof}

Now, let $l=k-1$ and by use Theorem \ref{thgen1}, it is easy to obtain 
the following.

\begin{corollary}
\label{cogen1}
For $k\geq 3$,
$$F_{1\mn2\mn\dots\mn (k-2)\mn(k-1)k}(x,y)=\sum_{j=0}^{k-1} (xF_{1\mn2\mn\dots\mn (k-2)\mn (k-1)k}(x,y))^j.$$
\end{corollary}

\begin{remark}
Similarly, the argument of $d$-increasing canonical decomposition, or 
the argumnet $d$-decreasing canonical decomposition yields others formulas, for example, 
formula for $F_{12\mn3\mn45}(x,y)$.
\end{remark}
%=========================================================================================
\section{Three letters pattern without internal dashes}
In this section, we give a complete answer for $F_\tau(x,y)$ where $\tau$ is 
a generalized pattern without internal dashes; that is, $\tau$ is $123$, $213$, 
$231$, $312$, and $321$. This by the following four subsections.

%************************
\subsection*{3.1 Patterns $123$ and $321$}

\begin{theorem}
$$F_{123}(x,y)=F_{321}(x,y)=\frac{1+xy-x+\sqrt{1-2x-3x^2-xy(2-2x-xy)}}{2x(x+y-xy)}.$$
\end{theorem}
\begin{proof}
Theorem \ref{th21c} yields, $F_{123}(x,y)=F_{321}(x,y)=H$ where 
	$$H=1+xH+\frac{x^2H^2}{1-xyH},$$
and the rest it is easy to see.
\end{proof}
%************************
\subsection*{3.2 Pattern $231$}

\begin{theorem}
$$F_{231}(x,y)=\frac{1-2x+2xy-\sqrt{1-4x+4x^2-4x^2y}}{2xy},$$
that is, for all $r,n\geq 0$
$$F_{231;r}(x)=\frac{1}{r+1}\binom{2r}{r} \frac{x^{2r+1}}{(1-2x)^{2r+1}},\quad 
  f_{231;r}(n)=\frac{2^{n-2r-1}}{r+1}\binom{n-1}{2r}\binom{2r}{r}.$$
\end{theorem}
\begin{proof}
Let $l=\pi^{-1}(n)$.
Since $\pi$ avoids $1\mn3\mn2$, each number in $\pi'$ is greater than
any of the numbers in $\pi''$. Therefore, $\pi'$ is a $1\mn3\mn2$-avoiding
permutation of the numbers $\{n-l+1,n-l+2,\dots,n-1\}$, while $\pi''$
is a $1-3-2$-avoiding permutation of the numbers $\{1,2,\dots,n-l\}$. On the
other hand, if $\pi'$ is an arbitrary $1\mn3\mn2$-avoiding
permutation of the numbers $\{n-l+1,n-l+2,\dots,n-1\}$ and $\pi''$ is an
arbitrary $1-3-2$-avoiding permutation of the numbers $\{1,2,\dots,n-l\}$, then
$\pi=(\pi',n,\pi'')$ is $1\mn3\mn2$-avoiding. \\

Now let us observe all posiiblities of $\pi'$ and $\pi''$ are empty or not. This 
yields 
	$$F_{231}(x,y)=1+x+2x(F_{231}(x,y)-1)+xy(F_{231}(x,y)-1)^2,$$
hence, the rest it is east to see.
\end{proof}
%********************************
\subsection*{3.3 Pattern $213$}

\begin{theorem}
$$F_{213}(x,y)=\frac{1-x^2+x^2y-\sqrt{1+2x^2-2x^2y+x^4-2x^4y+x^4y^2-4x}}{2x(1+xy-x)}.$$
\end{theorem}
\begin{proof}
Let $D(x,y)$ be the generating function of all $1\mn3\mn2$-avoiding permutations $(\alpha',n)\in S_n$ 
such that containing $213$ exactly $r$ times. Let $\alpha=(\alpha',n,\alpha'')$; if we consider the 
two cases $\alpha'$ empty or not we obtain $F_{213}(x,y)=1+D(x,y)F_{213}(x,y)$. Now let 
$\alpha=(\alpha',n)$; if we observe the two cases $\alpha'$ empty or not, then similarly 
it is easy to see 
	$$D(x,y)=x+x^2+x^2y(F_{213}(x,y)-1)+x^2(D(x,y)-1)+x^2(D(x,y)-1)(F_{213}(x,y)-1).$$
However,  
	$$F_{213}(x,y)=1+xF_{213}(x,y)\frac{1+x-xy+x(y-1)F_{213}(x,y)}{1-xF_{213}(x,y)},$$
hence, the rest it is easy to see.
\end{proof}
%*************************************
\subsection*{3.4 Pattern $312$}
\begin{theorem}
$$F_{312}(x,y)=\frac{1-x^2+x^2y-\sqrt{1+2x^2-2x^2y+x^4-2x^4y+x^4y^2-4x}}{2x(1+xy-x)}.$$
\end{theorem}
\begin{proof}
Let $\alpha\in\sss$; if $\alpha=\emptyset$, then there one permutation, otherwise 
by Proposition \ref{progen} we can write $\alpha=(\alpha^1,n,a_1,a_2,\dots,a_d,\alpha^{d+1},\alpha^d,\dots,\alpha^2)$ 
where all the entries of $\alpha^j$ greater than all the entries of $\alpha^{j+1}$, and 
$n>a_1>a_2>\dots>a_d$. Hence, for any $d=0,1$ the generating function of these permutations 
in these cases is $x^{d+1}F_{312}(x,y)$. Let $d\geq 2$, If $\alpha^{d+1}=\emptyset$, then 
the generating function of these permutations in this case is $x^{d+1}F_{312}^d(x,y)$, 
otherwise the generating function is $x^{d+1}yF_{312}^d(x,y)(F_{312}(x,y)-1)$. Hence 
$$F_{312}(x,y)=1+(x+x^2)F_{312}(x,y)+\sum_{d\geq2} x^{d+1}F_{312}^d(x,y)+\sum_{d\geq 2} x^{d+1}yF_{312}^d(x,y)(F_{312}(x,y)-1),$$
which means that
$$F_{312}(x,y)=1+xF_{312}(x,y)+\frac{x^2F_{312}(x,y)}{1-xF_{312}(x,y)}+\frac{x^2yF_{312}(x,y)(F_{312}(x,y)-1)}{1-xF_{312}(x,y)},$$
so the rest it is easy to see.
\end{proof}
%=========================================================================================
\section{Some of $3$-letter patterns with one dash}
In this section, we present some of examples $F_\tau(x,y)$  where $\tau$ is 
a generalized pattern with one dash. 
Theorem \ref{th21} yields as follows.

\begin{theorem} 
The generating function $F_{12\mn3}(x,y)$ 
is given by the continued fraction
$$
\frac1{1-\dfrac{x}{1-x+xy-\dfrac{xy}
{1-x+xy^2-\dfrac{xy^2}{\ddots}}}}.
$$
\end{theorem}

Theorem \ref{th21b} yields as follows.

\begin{theorem}
For any $k\geq 2$,
$$
F_{21\mn3}(x,y)=
1-\frac{x}
{x-\dfrac{1}
{1-\dfrac{x}
{xy-\dfrac{1}
{1-\dfrac{x}
{xy^2-\dfrac{1}{\ddots}}}}}}.$$
\end{theorem}

For $k=3$ and $l=2$ Theorem \ref{cogen1} yields the following.
\begin{theorem}
\label{th123}
$$F_{1\mn23}(x,y)=1+xF_{1\mn23}(x,y)+\sum_{d\geq 1} x^{d+1}y^{\binom{d}{2}}F_{1\mn23}(x,y)\prod_{j=0}^{d-1} F_{1\mn23}(xy^j,y).$$
\end{theorem}

\begin{corollary}
$$\begin{array}{ll}
F_{1\mn23;0}(x)&=\frac{1-x-\sqrt{1-2x-3x^2}}{2x^2};\\
F_{1\mn23;1}(x)&=\frac{x-1}{2x}+\frac{1-2x-x^2}{2x\sqrt{1-2x-3x^2}};\\
F_{1\mn23;2}(x)&=\frac{x^4}{(1-2x-3x^2)^{3/2}};\\
F_{1\mn23;3}(x)&=x^2-1+\frac{11x^7+43x^6+41x^5-7x^4-25x^3+x^2+5x-1}{(1-2x-3x^2)^{5/2}}.
\end{array}$$
\end{corollary}
\begin{proof}
By Theorem \ref{th123} with use $F_{1-23}(x,0)=F_{1-23;0}(x)$ we get 
	$$F_{1\mn23;0}(x)=1+xF_{1\mn23;0}(x)+x^2F_{1\mn23;0}^2(x),$$
which means the first formula holds.\\

By Theorem \ref{th123} we get 
	$$\frac{d}{dy} F_{1\mn23}(x,0)=x\frac{d}{dy} F_{1\mn23}(x,0)+2x^2F_{1\mn23}(x,0) \frac{d}{dy} F_{1\mn23}(x,0)+x^3F_{1\mn23}(x,0)^2F_{1\mn23}(0,0),$$
and with use $F_{1\mn23;1}(x)=\frac{d}{dy} F_{1\mn23}(x,y)\big{|}_{y=0}$ and the first formula, 
the second formula holds.\\

Similarly, by Theorem \ref{th123} and use $F_{1\mn23;r}(x)=\frac{1}{r!} \frac{d^r}{dy^r} F_{1\mn23}(x,y)\big{|}_{y=0}$ 
the others formulas holds.
\end{proof}

\begin{theorem}
\label{th213}
$$F_{2\mn13}(x,y)=\dfrac{1}{1-\dfrac{x}{1-\dfrac{x}{1-\dfrac{xy}{1-\dfrac{xy}{1-\dfrac{xy^2}{1-\dfrac{xy^2}{1-\ddots}}}}}}}.$$
\end{theorem}
\begin{proof}
By use Proposition \ref{progen} and Proposition \ref{progen1}, we obtain 
	$$F_{2\mn13}(x,y)=1+xF_{2\mn13}(x,y)\sum_{d\geq 0} x^dF^d_{2\mn13}(xy,y),$$
and the rest is easy to see.
\end{proof}
%==================================================================================
\section{Fruther results}
First of all, let us denote by $G_{\tau;\phi}(x,y)$ the generating function for 
the number of permutations in $S_n(1\mn3\mn2,\tau)$ such containing $\phi$ exactly $r$ times; 
that is 
$$G_{\tau;\phi}(x,z)=\sum_{n\geq 0} x^n\sum_{\pi\in S_n(1\mn3\mn2,\tau)} y^{a_\phi(\pi)}, $$ 
where $a_\phi(\pi)$ is the number of occurrences of $\phi$ in $\pi$. 
In this section, (similarly to above sections) we find 
$G_{\tau;\phi}(x,y)$ in terms of continued fractions or 
by explicit formula, for some cases of $\tau$ and $\phi$.

\begin{theorem}
The generating functions $G_{123;213}(x,y)$ and $G_{321;312}(x,y)$ are given by 
$$\dfrac{1}
{1-x-x^2(1-y)-\dfrac{x^2y}
{1-x-x^2(1-y)-\dfrac{x^2y}
{1-x-x^2(1-y)-\dfrac{x^2y}{\ddots}}}},$$
or, 
$$\frac{1-x-x^2+x^2y-\sqrt{(1-x-x^2)^2-2yx^2(1+x+x^2)+x^4y^2}}{2x^2y}.$$
\end{theorem}

\begin{theorem}
$$G_{123;231}(x,y)=H(x,y)+x^2(1-y)H(x,y)^2,$$
where $H(x,y)=\frac{1}{1-x-x^2yH(x,y)}$,  which means number of permutations in $S_n(1\mn3\mn2,123)$ 
such contains $231$ exactly $r\geq 1$ times is given by
	$$C_r\left(\binom{n}{2r}+\binom{n-1}{2r+1}\right) +C_{r-1}\binom{n-2}{2r-1},$$
where $C_r$ is the $r$th Catalan number, and for $r=0$ is given by $n$.
\end{theorem}

\begin{theorem}
The generating functions $G_{213;123}(x,y)$ and $G_{312;321}(x,y)$ given by 
	$$\frac{1-x-x^2+xy-\sqrt{(1-x-x^2)^2-2xy(1-x+x^2)+x^2y^2}}{2xy(1-x)}.$$
\end{theorem}

As a concluding remark we note that are many questions left to answer. If there 
exist a bijection between, for example, the set of $1\mn3\mn2$-avoiding permutations 
in $S_n$ such containing $r$ times the generalized pattern $21\mn3\mn\dots\mn k$, and  
the set of $1\mn3\mn2$-avoiding permutations in $S_n$ such that containing $r$ times 
the classically pattern $1\mn2\mn3\mn\dots\mn k$, where $r=0,1,\dots,k-2$.
%======================================================================

\end{document}